\documentclass[reqno,12pt]{amsart}

\newtheorem{rem}{Remark}

\newtheorem{prop}{Proposition}

\newtheorem{df}{Definition}

\newcommand\eps\varepsilon
\newcommand\ph\varphi
\newcommand\kap\Lambda

\usepackage{enumerate}
\usepackage{graphicx}
\usepackage{float}

\begin{document}

\title[On Periodic Solutions to Lagrangian System]
{ On Periodic Solutions to Lagrangian System With Singularities }

\author[Oleg Zubelevich]{Oleg Zubelevich\\ \\\tt
 Dept. of Theoretical mechanics,  \\
Mechanics and Mathematics Faculty,\\
M. V. Lomonosov Moscow State University\\
Russia, 119899, Moscow,  MGU \\
 }
\date{}
\thanks{Partially supported by grant
 RFBR  18-01-00887}

\begin{abstract}A  Lagrangian system  with singularities  is considered. The configuration space  is a non-compact manifold that depends on time. A set of  periodic solutions has been found.
\end{abstract}

\maketitle
\numberwithin{equation}{section}
\newtheorem{theorem}{Theorem}[section]
\newtheorem{lemma}[theorem]{Lemma}
\newtheorem{definition}{Definition}[section]

\section{Introduction}
Let us start from a model example.

Consider the plane  $\mathbb{R}^2$ with the standard Cartesian frame $$Oxy,\quad \boldsymbol r=x\boldsymbol e_x+y\boldsymbol e_y$$ and the standard  Euclidean norm $|\cdot|$.

A particle of mass $m$ moves in the plane being influenced by a force with the potential 
$$V(\boldsymbol r)=-\gamma\Big(\frac{1}{|\boldsymbol r-\boldsymbol r_0|^n}+\frac{1}{|\boldsymbol r+\boldsymbol r_0|^n}\Big),$$ 
here $\gamma>0$ is a constant; $\boldsymbol r_0\ne 0$ is a given vector.
The Lagrangian of this system is as follows
\begin{equation}\label{dtgt66}L(\boldsymbol r,\boldsymbol {\dot r})=\frac{m}{2}|\boldsymbol {\dot r}|^2-V(\boldsymbol r).\end{equation}

If $n=1$ then this system is  classical  the two fixed centers problem and it is integrable \cite{arn}. The case of potentials of type  $V(\boldsymbol r)\sim -1/|\boldsymbol r|^2$ is mainly used in  quantum mechanics but the classical statement is also studied, see \cite{Martinez}, \cite{Krishnaswami} and references therein. 

One of the consequences from the main Theorem \ref{main_theor} (see  the next section) is as follows.

\begin{prop}\label{xvfvvvv} For any $\omega>0,\quad n\ge 2$ and for any $m\in\mathbb{N}$  system (\ref{dtgt66}) has an $\omega$-periodic solution $\boldsymbol r(t)$ 
that  first $m$ times coils   clockwise  about the point $\boldsymbol r_0$  and then $m$ times coils  counterclockwise about the point $-\boldsymbol r_0$ . If $m=1$ then the solution $\boldsymbol r(t)$  forms an $8$-like curve in the plane. \end{prop}

This fact does not hold for $n=1$ \cite{ger}.

Consider another pure classical mechanics example. Introduce in our space an inertial  Cartesian frame $Oxyz$ such that the gravity is $\boldsymbol g=-g\boldsymbol e_z$. 

Let a particle of mass $m$ be influenced by gravity and  slides without friction  on a surface 
$$z=f(\boldsymbol r)=-\gamma\Big(\frac{1}{|\boldsymbol r-\boldsymbol r_0|^n}+\frac{1}{|\boldsymbol r+\boldsymbol r_0|^n}\Big),\quad \boldsymbol r=x\boldsymbol e_x+y\boldsymbol e_y.$$
The corresponding Lagrangian is
\begin{equation}\label{fgrgtcc}L(\boldsymbol r,\boldsymbol {\dot r})=\frac{m}{2}\Big(|\boldsymbol{\dot r}|^2+\big(\nabla f(\boldsymbol r),\boldsymbol{\dot r}\big)^2\Big)-mgf(\boldsymbol r).\end{equation}

Proposition \ref{xvfvvvv} holds for system (\ref{fgrgtcc}) also.

Existence problems for periodic solutions to Lagrangian systems have intensively been studied since  the beginning of the 20th century  and even earlier. There is an immense  number of different results and methods developed in this field. We mention only a few of them which are mach closely related to this article.

In \cite{capozzi} periodic solutions have been obtained for the Lagrangian system with Lagrangian
$$L(t,x,\dot x)=\frac{1}{2}g_{ij}(x)\dot x^i\dot x^j-W(t,x),\quad x=(x^1,\ldots,x^m)\in\mathbb{R}^m$$
here and in the sequel we use the Einstein summation convention. The form $g_{ij}$ is symmetric and positive definite: 
$$g_{ij}\xi^i\xi^j\ge \mathrm{const}_1\cdot|\xi|^2.$$ 
The potential is as follows $W(t,x)=V(x)+g(t)\sum_{i=1}^m x^i,$ where $V$ is a bounded function $|V(x)|\le\mathrm{const}_2$ and $g$ is an $\omega-$ periodic function. The functions $V,g_{ij}$ are even. 

Under these assumptions the authors prove that there exists a nontrivial $\omega-$ periodic solution.

Our main tool to obtain periodic solutions is variation technique.
Variational problems and Hamiltonian systems have been studied extensively. Classic references of these
subjects are \cite{23}, \cite{28}, \cite{14}.

\section{The Main Theorem}
Introduce some notations.
Let $x=(x^1,\ldots, x^m)$ and $\ph=(\ph^1,\ldots,\ph^n)$ be points of the standard $\mathbb{R}^m$ and $\mathbb{R}^n$ respectively.
Then let $z$ stand for the point $(x,\varphi)\in \mathbb{R}^{m+n}$. By $|\cdot|$ denote the standard Euclidean  norm of $\mathbb{R}^k,\quad k=m,m+n,$ that is $|x|^2=\sum_{i=1}^k(x^i)^2.$

The variable $z$ can consist only of $x$ without $\varphi$ or conversely. 

Assume that we are provided with a set $\sigma\in\mathbb{R}^{m+n}$ such that

1) If $z=(x,\varphi)\in\sigma$ then $-z\in\sigma$ and $(x,\varphi+2\pi p)\in\sigma$ for any 
$p=(p^1,\ldots ,p^n)\in\mathbb{Z}^n$;  

2) the set $\sigma$ does not have accumulation points.

Introduce a manifold $\Sigma=\mathbb{R}^{m+n}\backslash\sigma.$

\begin{rem}\label{remqq}Actually, the space $\mathbb{R}^{m+n}$ wraps  a cylinder $\mathcal C=\mathbb{R}^m\times \mathbb{T}^n$, where $\mathbb{T}^n$ is the torus with angle variables $\varphi$. By the same reason, $\Sigma$ can be regarded as the cylinder with several dropped away points.
Nevertheless for analysis purposes we prefer to deal with $\mathbb{R}^{m+n}$. \end{rem}

The main object of our study is the following Lagrangian system with Lagrangian
\begin{equation}\label{dtgd5}L(t,z,\dot z)=\frac{1}{2}g_{ij}\dot z^i\dot z^j+a_i\dot z^i-V,\quad z=(z^1,\ldots,z^{m+n}).\end{equation}
\begin{rem}The term $a_i\dot z^i$ in the Lagrangian corresponds to the so called gyroscopic forces.  For example, the Coriolis force and the Lorentz force are gyroscopic. 
\end{rem}

 The functions $g_{ij},a_i$ depend on $(t,z)$ and belong to $C^2(\mathbb{R}\times \Sigma)$; moreover all these functions are $2\pi-$periodic in each variable $\varphi^j$ and $\omega-$periodic in the variable $t,\quad \omega>0$. For all $(t,z)\in \mathbb{R}\times \Sigma$ it follows that  $g_{ij}=g_{ji}$.
 
 The function $V(t,z)\in C^2(\mathbb{R}\times\Sigma)$  is also $2\pi-$periodic in each variable $\varphi^j$ and $\omega-$periodic in the variable $t$. 

We also assume that there are positive constants $C,M,A,K,P$ such that for all $(t,z)\in\mathbb{R}\times \Sigma$ and $ \xi\in \mathbb{R}^{m+n}$  we have
\begin{equation}\label{sdfff}|a_i(t,z)|\le C+M|z|,\quad \frac{1}{2}g_{ij}(t,z)\xi^i\xi^j\ge K|\xi|^2;\end{equation}
for all $(t,z)\in\mathbb{R}\times\Sigma,\quad \sigma'\in\sigma$ it follows that
\begin{equation}\label{q_sdfff}V(t,z)\le A|z|^2-\frac{P}{|z-\sigma'|^2}+C_1.\end{equation}

The Lagrangian $L$ is defined up to an additional constant (see Definition \ref{dfg5} below), so the constant $C_1$ is not essential.

System (\ref{dtgd5}) obeys  the following ideal constraints:
\begin{equation}\label{constr}
f_j(t,z)=0,\quad j=1,\ldots,l<m+n,\quad f_j\in C^2(\mathbb{R}^{m+n+1}).\end{equation}
The functions $f_j$ are also $2\pi-$periodic in each variable $\varphi^j$ and $\omega-$periodic in the variable $t$.

Introduce a set $$F(t)=\{z\in\mathbb{R}^{m+n}\mid f_j(t,z)=0,\quad j=1,\ldots,l\}.$$
 Assume that 
\begin{equation}\label{dfggfgfg3e}
\mathrm{rank}\,\Big(\frac{\partial f_j}{\partial z^k}(t,z)\Big)=l\end{equation}for all $z\in F(t)$. So that $F(t)$ is a smooth manifold.

Assume also that all the functions $f_j$ are either  odd:
\begin{equation}\label{ncdnvv1}
f_j(-t,-z)=-f_j(t,z)\end{equation}
or even 
\begin{equation}\label{ncdnvv2}
f_j(-t,-z)=f_j(t,z).\end{equation}

The set $\sigma$ belongs to $F(t)$:
\begin{equation}\label{dvt5y56}\sigma\subset F(t).\end{equation}

\begin{rem} Actually, it is sufficient to say that all the functions   are defined and have formulated above properties in some open symmetric vicinity of the manifold $F(t)$.
We believe that this generalization is unimportant and keep referring to the whole space $\mathbb{R}^{m+n}$ just for simplicity of exposition. 
\end{rem}
\begin{rem}
The set $\sigma$ can be void. In this case the condition (\ref{dvt5y56}) is  dropped and formula (\ref{q_sdfff}) takes the form
$$V(t,z)\le A|z|^2+C_1.$$

\end{rem}
\begin{df}\label{dfg5}We shall say that a function $z(t)\in C^2(\mathbb{R},\Sigma)$ is a solution to system (\ref{dtgd5}), (\ref{constr}) if there exists a set of functions $$\{\alpha^1,\ldots,\alpha^l\}\subset C(\mathbb{R})$$ such that
\begin{equation}\label{xvvv}\frac{d}{dt}\frac{\partial L}{\partial \dot z^i}(t,z(t),\dot z(t))-\frac{\partial L}{\partial  z^i}(t,z(t),\dot z(t))=\alpha^j(t)\frac{\partial f_j}{\partial z^i}(t,z(t))\end{equation}
and \begin{equation}\label{cfgbgg}z(t)\in F(t)\end{equation}
for all real $t$.

In the absence of constraints (\ref{constr}) right side of (\ref{xvvv}) is equal to zero and the condition (\ref{cfgbgg}) is dropped.
\end{df}The functions $\alpha^i$  are defined from equations (\ref{xvvv}), (\ref{constr}) uniquely.
\begin{df}Let $$\mathcal V_p\subset C(\mathbb{R},\Sigma),\quad p=(p^1,\ldots,p^n)\in\mathbb{Z}^n$$ stand for a set of functions $z(t)=(x(t),\varphi(t))$ that satisfy 
the following properties for all $t\in\mathbb{R}$:

1) $x(t+\omega)=x(t);\quad \varphi(t+\omega)=\varphi(t)+2\pi p,$

2) $z(t)\in F(t)$.

In the absence of constraints (\ref{constr})  condition 2) is omitted.
\end{df}

\begin{df} 
We shall say that two functions $z_1,z_2\in\mathcal V_p$ are homotopic to each other  iff there exists a mapping
$z(s,t)\in C([0,1]\times\mathbb{R},\Sigma)$ such that

1) $z(s,\cdot)\in\mathcal V_p,$

2) $z(0,t)=z_1(t),\quad z(1,t)=z_2(t).$

By $[z]$ we denote the homotopy class of the function $z$.
\end{df}

\begin{theorem}\label{main_theor}
Assume that 

1) all the functions are even:
$$g_{ij}(-t,-z)=g_{ij}(t,z),\quad a_i(-t,-z)=a_i(t,z),\quad V(-t,-z)=V(t,z);$$

2) the following inequality holds
\begin{equation}\label{xcxddfd}K-\frac{M\omega}{\sqrt 2}-\frac{A\omega^2}{2}>0;\end{equation}

3) for some $\nu\in\mathbb{Z}^n$ the set $\mathcal V_\nu$ is non-void and  there is a function
$$ \tilde z(t)\in\mathcal V_\nu\cap C^1(\mathbb{R},\Sigma);\quad \tilde z(-t)=-\tilde z(t);$$

4) and
\begin{equation}\label{cvbvggdfgkk}
\inf\{\|z-\sigma'\|_{C[0,\omega]}\mid z\in [\tilde z],\quad \sigma'\in\sigma\}>0.\end{equation}

Then 
system (\ref{dtgd5}), (\ref{constr}) has an odd solution $$z_*(t)\in [\tilde z],\quad z_*(-t)=-z_*(t).$$

This assertion remains valid in the absence of constraints (\ref{constr}).  When $\sigma=\emptyset$ then the assertion remains valid with condition 4) omitted. 
\end{theorem}
Actually the solution stated in this theorem is as smooth as it is allowed by smoothness of the Lagrangian $L$ and the functions $f_j$ up to $C^\infty$.

Loosely speaking,  condition 4) of this theorem (formula (\ref{cvbvggdfgkk})) implies that we look for a solution among the curves those can not be shrunk into a point. See for example problem from the Introduction.
\begin{rem}
If all the functions do not depend on $t$ then we can choose $\omega$ to be arbitrary small and inequality (\ref{xcxddfd}) is satisfied.  Taking a vanishing sequence of $\omega$, we obtain infinitely many periodic solutions of the same homotopic type.
\end{rem}
\begin{rem}\label{dtrgrr}Condition 2) of the Theorem is essential. 
Indeed,  system 
$$L(t,x,\dot x)=\frac{1}{2}\dot x^2-\frac{1}{2}\Big(x-\sin t\Big)^2$$
obeys all the conditions except inequality (\ref{xcxddfd}). It is easy to see that the corresponding equation  $\ddot x+x=\sin t$ does not have periodic solutions.

\end{rem}

\subsection{Examples}
Our first example is as follows.
\begin{figure}[H]
\centering
\includegraphics[width=45mm, height=30 mm]{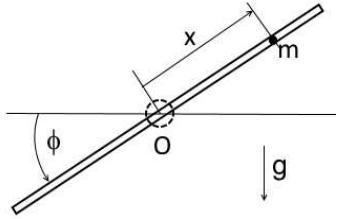}
\caption{the tube and the ball \label{overflow}}
\end{figure}

{\it A thin  tube can rotate freely in the vertical plane about a fixed horizontal axis passing through its centre $O$. A moment of inertia of the tube about this axis is equal to $J$. The mass of the tube is distributed symmetrically such that tube's centre of mass is placed at the point $O$.

Inside the tube there is a small ball  which can slide  without friction. The mass of the ball is $m$. The  ball can pass by the point $O$ and fall out from the ends of the tube.

The system  undergoes  the standard gravity field $\boldsymbol g$.}

It seems to be evident that for typical motion the ball reaches an end of the tube and falls down out the tube. It is surprisingly, at least for the first glance, that this system has very many periodic solutions such that the tube turns around several times during the period. 

The sense of generalized coordinates $\phi,x$ is clear from  Figure \ref{overflow}.

The Lagrangian of this system is as follows
\begin{equation}\label{Lagr}L(x,\phi,\dot x,\dot \phi)= \frac{1}{2}\Big(mx^2+J\Big)\dot\phi^2+\frac{1}{2}m\dot x^2-mgx\sin\phi.\end{equation}
Form theorem \ref{main_theor} it follows that 
for any constant $\omega>0$ system (\ref{Lagr}) has a solution $\phi(t),x(t),\quad t\in\mathbb{R}$ such that 

1) $x(t)=-x(-t),\quad  \phi(t)=-\phi(-t);$

2) $x(t+\omega)=x(t),\quad \phi(t+\omega)=\phi(t)+2\pi .$

This result shows that for any $\omega>0$  the system has an $\omega-$periodic motion such that the tube turns around  once during the period. The length of the tube should be chosen properly. 

Our second example is a counterexample. Let us show that the first condition of the theorem \ref {main_theor}  can not be omitted. 

{\it Consider a mass point $m$ that slides on a right circular cylinder  of radius $r$. The surface of the cylinder is perfectly smooth. The axis $x$ of the cylinder is parallel to the gravity $\boldsymbol g$ and directed upwards.}

The Lagrangian of this system is 
\begin{equation}\label{xvfvffvvf}L(x,\varphi,\dot x,\dot\varphi)=\frac{m}{2}\Big(r^2\dot \varphi^2+\dot x^2\Big)-mgx.\end{equation}
All the conditions except the evenness  are satisfied but it is clear this system does not have periodic solutions.

\section{ Proof of Theorem \ref{main_theor} }
In this section we use several standard facts from functional analysis and the Sobolev spaces theory \cite{edv}, \cite{adams}.

By $c_1,c_2,\ldots$ we denote inessential positive constants.
\subsection{}
Recall that the Sobolev space $H^1_{\mathrm{loc}}(\mathbb{R})$ consists of  functions $u(t),t\in\mathbb{R}$ such that $u,\dot u\in L^2_{\mathrm{loc}}(\mathbb{R})$. The following embedding holds $H^1_{\mathrm{loc}}(\mathbb{R})\subset C(\mathbb{R})$.

Recall another standard fact.
\begin{lemma}\label{xddddd}
Let $u\in H^1_{\mathrm{loc}}(\mathbb{R})$ and $u(0)=0$. Then for any $a>0$ we have
$$\|u\|^2_{L^2(0,a)}\le \frac{a^2}{2}\|\dot u\|^2_{L^2(0,a)},\quad \|u\|^2_{C[0,a]}\le a \|\dot u\|^2_{L^2(0,a)}.$$\end{lemma}
Here and below the notation $\|\dot u\|_{L^2(0,a)}$ implies that 
$$\Big\|\dot u\Big|_{(0,a)}\Big\|_{L^2(0,a)}$$ the same is concerned to $\|u\|_{C[0,a]}$ etc.

\subsection{}Here we collect several spaces which are needed in the sequel.
\begin{df}
By $X$ denote a space of functions $u\in   H^1_{\mathrm{loc}}(\mathbb{R})$ such that for all $t\in\mathbb{R}$ the following conditions hold
$$u(-t)=-u(t),\quad u(t+\omega)=u(t).$$
\end{df}
By virtue of lemma \ref{xddddd}, the mapping $u\mapsto \|\dot u\|_{L^2(0,\omega)}$ determines a norm in $X$. This norm is denoted by $\|u\|$. The norm $\|\cdot\|$ is equivalent to the standard norm of $H^1[0,\omega]$. The space $(X,\|\cdot\|)$ is a Banach space. Since the norm $\|\cdot\|$ is generated by an inner product
$$(u,v)_X=\int_0^{\omega}\dot u(t)\dot v(t)dt$$ the space $X$ is also a real Hilbert space, particularly this implies that $X$ is a reflexive Banach space.

\begin{df} Let $\Phi$ stand for the space $\{ct+u(t)\mid c\in\mathbb{R},\quad u\in X\}$. \end{df}
By the same argument, $(\Phi,\|\cdot\|)$ is a reflexive Banach space. Observe also that $\Phi=\mathbb{R}\oplus X$ and by direct calculation 
we get
$$\|\psi\|^2=\omega c^2+\|u\|^2,\quad \psi(t)=ct+u(t)\in\Phi.$$
Observe that $X\subset \Phi$.
\begin{df}
Let $E$ stand for the space
 $$X^m\times\Phi^n=\{z(t)=(x^1,\ldots, x^m,\varphi^1,\ldots,\ph^n)(t)\mid x^i\in X,\quad \varphi^j\in\Phi\}.$$\end{df}
 The space $E$ is also a real Hilbert space
with an  inner product defined as follows
$$(z,y)_E=\int_0^{\omega}\sum_{i=1}^{m+n}\dot z^i(t)\dot y^i(t)dt,$$
where $ z=(z^k),y=(y^k)\in E,\quad k=1,\ldots, m+n.$

We denote the corresponding norm  in $E$  by the same symbol and write
$$\|z\|^2=\big\||z|\big\|^2=\sum_{k=1}^{m+n}\|z^k\|^2.$$
The space $E$ is also a reflexive Banach space.

Introduce the following set 
$$E_0=\Big\{(x,\varphi)\in E\Big|\varphi^j=\frac{2\pi \nu_j}{\omega}t+u_j,\quad \forall x\in X,\quad \forall u_j\in X,\quad j=1,\ldots,n\Big\}.$$
This set is a closed plane of codimension $n$ in $E$.

If $(x,\varphi)\in E_0$ then $\varphi(t+\omega)=\varphi(t)+2\pi \nu.$

\begin{df} Let $Y$ stand for the space $$\big\{u\in L^2_{\mathrm{loc}}(\mathbb{R})\mid u(t)=u(-t),\quad u(t+\omega)=u(t)\quad \mbox{almost everywhere in }\mathbb{R}\big\}.$$ \end{df}
The space $Y^{m+n}$ is a Hilbert space with respect to the inner product
$$(z,y)_Y=\int_0^{\omega}\sum_{i=1}^{m+n} z^i(t) y^i(t)dt,$$
where $ z=(z^k),y=(y^k)\in Y^{m+n},\quad k=1,\ldots, m+n.$

\subsection{} Fix a positive constant $K_1$ and consider a function 
\begin{equation}\label{vfdgdfg66}u:[0,\omega]\to\Sigma,\quad u(t)=(u^1,\ldots,u^{m+n})(t),\quad u^i\in H^1[0,\omega]\end{equation} such that
the following two conditions hold
\begin{equation}\label{fgggfdvv}
\int_0^\omega\frac{dt}{|u(t)-\hat u|^2}\le K_1,\quad \|u\|_{H^1[0,\omega]}\le K_1\end{equation}
with some constant vector $\hat u\in \mathbb{R}^{m+n}$.
\begin{lemma}\label{xvcvcccc}
There exists a positive constant $K_2$ such that for any function $u$ that satisfies (\ref{vfdgdfg66}) and (\ref{fgggfdvv}) one has
$$|u(t)-\hat u|\ge K_2\sqrt{\sum_{k=1}^{n+m}\|u^k-\hat u^k\|^2_{C^{0,1/2}[0,\omega]}},\quad t\in[0,\omega],\quad j=1,\ldots,N.$$The constant $K_2$ depends only on $K_1$. \end{lemma}Here $C^{0,1/2}[0,\omega]$ is the Holder space.

{\it Proof of Lemma \ref{xvcvcccc}.} By $c_1,c_2,\ldots$ we denote inessential positive constants.

Observe  that $H^1[0,\omega]\subset C^{0,1/2}[0,\omega]$ and \begin{equation}\label{xfvbvvv}
\|\cdot\|_{C^{0,1/2}[0,\omega]}\le c_1\|\cdot\|_{H^1[0,\omega]}.\end{equation}

Take arbitrary $t'\in[0,\omega]$.  Fix $j$ and introduce the following notations $v=(v^1,\ldots, v^{m+n}),$
$$  v^k=u^k-\hat u^k,\quad a=|v(t')|,\quad   b=\sqrt{\sum_{k=1}^{n+m}\|v^k\|^2_{C^{0,1/2}[0,\omega]}}.$$
Observe that from conditions of the Lemma and formula (\ref{xfvbvvv}) 
we get\begin{equation}\label{vfvfgffffhdfgh}b\le c_3.\end{equation}
It follows that
$$|v^k(t)-v^k(t')|\le \|v^k\|_{C^{0,1/2}[0,\omega]}|t-t'|^{1/2},\quad t\in[0,\omega]$$ and
$$|v^k(t)|\le |v^k(t')|+\|v^k\|_{C^{0,1/2}[0,\omega]}|t-t'|^{1/2}.$$Consequently, we obtain
$$|v(t)|\le a+b|t-t'|^{1/2}.$$
We then have
$$
K_1\ge\int_0^\omega\frac{dt}{|v(t)|^2}\ge \int_0^{t'}\frac{dt}{\big(a+b(t'-t)^{1/2}\big)^2}+
\int_{t'}^\omega\frac{dt}{\big(a+b(t-t')^{1/2}\big)^2}.$$
The last two integrals are computed easily with the help of the change $\xi=a+b(t'-t)^{1/2}$ for the first integral, and 
$\xi=a+b(t-t')^{1/2}$ for the second integral respectively.

Assume for definiteness that $t'\in[0,\omega/2]$. Other case is carried out analogously. Then inequality 
$$ \int_{t'}^\omega\frac{dt}{\big(a+b(t-t')^{1/2}\big)^2}\le K_1$$ and formula (\ref{vfvfgffffhdfgh}) imply
$$a\ge c_2b.$$

The Lemma is proved.

\begin{df}Let $W$ stand for a set $E_0\cap \mathcal [\tilde z].$\end{df}
\subsection{}

Introduce the Action Functional $S:W\to\mathbb{R},$
$$S(z)=\int_0^{\omega}L(t,z,\dot z)dt.$$
Our next goal is to prove that this functional attains its minimum.

By using estimates (\ref{sdfff}), (\ref{q_sdfff}) we get
$$S(z)\ge \int_0^{\omega}\big(K|\dot z|^2-|\dot z|(C+M|z|)-A|z|^2\big)dt.$$
From the Cauchy inequality and Lemma \ref{xddddd} it follows that
\begin{align}
 \int_0^{\omega}|\dot z|| z|dt&\le\frac{\omega}{\sqrt 2}\|z\|^2, \nonumber\\
\int_0^{\omega}| z|^2dt&\le \frac{\omega^2}{2}\|z\|^2, \nonumber\\
\int_0^{\omega}| \dot z|dt&\le \sqrt{\omega}\|z\|.\nonumber
\end{align}
We finally obtain 
\begin{equation}\label{xfffff}
S(z)\ge \Big(K-\frac{M\omega}{\sqrt 2}-\frac{A\omega^2}{2}\Big)\|z\|^2-C\sqrt{\omega}\|z\|.
\end{equation}

By formula (\ref{xcxddfd}) the functional $S$ is coercive:
\begin{equation}\label{dvfdffvfv}S(z)\to \infty\end{equation}as $\|z\|\to\infty.$

Note that the Action Functional which corresponds to system (\ref{xvfvffvvf}) is also coercive but, as we see above, property (\ref{dvfdffvfv}) by itself  does not imply existence results.

\subsection{} Let $\{z_k\}\subset W$ be a minimizing sequence:
$$S(z_k)\to \inf_{z\in W} S(z)$$ as $k\to\infty.$

By formula (\ref{dvfdffvfv}) the sequence $\{z_k\}$ is bounded: $\sup_k\|z_k\|<\infty.$ 
\begin{lemma}\label{qwwfdf}
The following estimates hold
\begin{align}
\inf\{|z_k(t)-&\sigma'|\mid t\in\mathbb{R},\quad \sigma'\in\sigma,\quad k\in\mathbb{N}\}\nonumber\\
&=\inf\{|z_k(t)-\sigma'|\mid t\in[0,\omega],\quad \sigma'\in\sigma,\quad k\in\mathbb{N}\}>0.\label{vgffdddd}\end{align}\end{lemma}
{\it Proof of lemma \ref{qwwfdf}.} The middle equality in formula (\ref{vgffdddd}) follows from the fact that $z_k\in\mathcal V_\nu$, we will prove inequality.

It is sufficient to check inequality (\ref{vgffdddd}) for $\sigma'\in \sigma\cap\Lambda,$
$$\Lambda=\bigcup_{k\in\mathbb{N}}\{z\in\mathbb{R}^{m+n}\mid |z-z_k(t)|\le 1,\quad t\in[0,\omega]\}.$$ The set $ \sigma\cap\Lambda$ is finite.

By formula (\ref{q_sdfff}) we obtain
\begin{align}
P\int_0^\omega&\frac{dt}{|z_k(t)-\sigma'|^2}+\Psi_k\nonumber\\
&\le c_7-\int_0^\omega V(t,z_k(t))dt+\Psi_k\le S(z_k)+c_7.\label{xvfffff}\end{align}
Here 
$$\Psi_k=\int_0^\omega a_i(t,z_k(t))\dot z_k(t)dt,\quad c_7=A\|z_k\|^2_{L^2[0,\omega]}+C_1\omega.$$ By formula (\ref{sdfff}) and lemma \ref{xddddd} this sequence is bounded:
$$|\Psi_k|\le c_5\|z_k\|+c_6\|z_k\|^2.$$

So that from formulas  (\ref{xvfffff}), (\ref{cvbvggdfgkk}) and  Lemma \ref{xvcvcccc}, Lemma \ref{qwwfdf} is proved.

Since the space $E$ is reflexive, this sequence contains a weakly convergent subsequence. Denote this subsequence in the same way: $z_k\to z_*$ weakly in $E$.

Moreover,  the space $H^1[0,\omega]$ is compactly embedded in $C[0,\omega]$. Thus extracting a subsequence from the subsequence and keeping the same notation we also have
\begin{equation}\label{scr 4vt}\max_{t\in[0,\omega]}|z_k(t)-z_*(t)|\to 0,\end{equation}
as $k\to\infty$;
and
$$\inf\{|z_*(t)-\sigma'|\mid t\in\mathbb{R},\quad \sigma'\in\sigma\}>0.$$

The set $E_0$ is convex and strongly closed therefore it is weakly closed: $z_*\in E_0$. By continuity (\ref{scr 4vt}) one also gets $z_*\in W$ (see Lemma \ref{xzsdffgff} in the Appendix).

\subsection{} Let us show that $\inf_{z\in W} S(z)=S(z_*).$
\begin{lemma}\label{xvvfvv}Let a sequence $\{u_k\}\subset \Phi$
weakly converge to $u\in\Phi$ (or $u_k,u\in X$ and $u_k\to u$ weakly in $X$); and also $\max_{t\in [0,\omega]}|u_k(t)- u(t)|\to 0$ as $k\to\infty$.

Then for any $f\in C(\mathbb{R})$ and for any $v\in L^2(0,\omega)$ it follows that
$$\int_0^{\omega}f(u_k)\dot u_k vdt\to \int_0^{\omega}f(u)\dot u vdt,\quad \mbox{as}\quad k\to\infty.$$
 \end{lemma}
 Indeed,
 \begin{align}
 \int_0^{\omega}&f(u_k)\dot u_k vdt\nonumber\\
&=\int_0^{\omega}\big(f(u_k)-f(u)\big)\dot u_k vdt+\int_0^{\omega}f(u)\dot u_k vdt.\nonumber\end{align}
The function $f$ is uniformly continuous in a compact set $$\Big[\min_{t\in[0,\omega]}\{u(t)\}-c,\max_{t\in[0,\omega]}\{u(t)\}+c\Big]$$ with some constant $c>0$.
Consequently we obtain $$\max_{t\in [0,\omega]}|f(u_k(t))- f(u(t))|\to 0.$$

Since the sequence $\{u_k\}$ is weakly convergent it is bounded: $$\sup_k\|u_k\|<\infty$$ particularly, we get
$$\|\dot u_k\|_{L^2(0,\omega)}<\infty.$$
So that as $k\to\infty$
\begin{align}
\Big|&\int_0^{\omega}\big(f(u_k)-f(u)\big)\dot u_k vdt\Big|\nonumber\\
&\le \big\|v\big(f(u_k)- f(u)\big)\big\|_{L^2(0,\omega)}\cdot \|\dot u_k\|_{L^2(0,\omega)}\to 0.\nonumber
\end{align}
To finish the proof it remains to observe that a function
$$w\mapsto\int_0^{\omega}f(u)\dot w vdt$$ belongs to $\Phi'$ (or to $X'$). Indeed,
$$\Big|\int_0^{\omega}f(u)\dot w vdt\Big|\le \max_{t\in[0,\omega ]}|f(u(t))|\cdot\|v\|_{L^2(0,\omega)}\|w\|.$$

\subsection{} The following lemma is proved similarly.
\begin{lemma}\label{xvvfvv1}Let a sequence $\{u_k\}\subset \Phi$ (or $\{u_k\}\subset X$)
be such that $$\max_{t\in [0,\omega]}|u_k(t)- u(t)|\to 0$$
as $k\to\infty.$

Then for any $f\in C(\mathbb{R})$ and for any $v\in L^1(0,\omega)$ it follows that
$$\int_0^{\omega}f(u_k) vdt\to \int_0^{\omega}f(u) vdt$$
as $k\to\infty.$
 \end{lemma}
\subsection{}
Introduce a function
$p_k(t,\xi)=L(t,z_k,\dot z_*+\xi)$. The function $p_k$ is a quadratic polynomial of $\xi\in\mathbb{R}^{m+n}$, so that
$$p_k(t,\xi)=L(t,z_k,\dot z_*)+\frac{\partial L}{\partial \dot z^i}(t,z_k,\dot z_*)\xi^i+\frac{1}{2}\frac{\partial^2 L}{\partial \dot z^j\partial \dot z^i}(t,z_k,\dot z_*)\xi^i\xi^j.$$
The last term in this formula is non-negative:
$$\frac{\partial^2 L}{\partial \dot z^j\partial \dot z^i}(t,z_k,\dot z_*)\xi^i\xi^j=g_{ij}(t,z_k)\xi^i\xi^j\ge 0.$$
We consequently obtain
$$p_k(t,\xi)\ge L(t,z_k,\dot z_*)+\frac{\partial L}{\partial \dot z^i}(t,z_k,\dot z_*)\xi^i.$$
It follows that
\begin{align}S(z_k)&=\int_0^{\omega}p_k(t,\dot z_k-\dot z_*)dt\ge \int_0^{\omega}L(t,z_k,\dot z_*)dt\nonumber\\&+
\int_0^{\omega}\frac{\partial L}{\partial \dot z^i}(t,z_k,\dot z_*)(\dot z_k^i-\dot z_*^i)dt.\label{xfbf}\end{align}
From Lemma \ref{xvvfvv} and  Lemma \ref{xvvfvv1} it follows that 
$$\int_0^{\omega}L(t,z_k,\dot z_*)dt\to \int_0^{\omega}L(t,z_*,\dot z_*)dt\quad \mbox{as}\quad k\to\infty,$$
and
$$\int_0^{\omega}\frac{\partial L}{\partial \dot z^i}(t,z_k,\dot z_*)(\dot z_k^i-\dot z_*^i)dt\to 0\quad \mbox{as}\quad k\to\infty.$$
Passing to the limit as $k\to\infty$ in  (\ref{xfbf}) we finally yield 
$$\inf_{z\in W}S(z)\ge S(z_*)\Longrightarrow \inf_{z\in W}S(z)= S(z_*).$$
\begin{rem}
Basing upon these formulas one can estimate the norm $\|z_*\|$. Indeed, take a function $\hat z\in W$ then due to formula (\ref{xfffff}) one obtains
$$S(\hat z)\ge S(z_*)\ge \Big(K-\frac{M\omega}{\sqrt 2}-\frac{A\omega^2}{2}\Big)\|z_*\|^2-C\sqrt{\omega}\|z_*\| ,$$
here $S(\hat z)$ is an explicitly calculable  number. \end{rem}

\subsection{}
Now from this point we begin proving the theorem under the assumption that the constraints are odd (\ref{ncdnvv1}).

Thus for any $v\in X^{m+n}$ such that
$$\frac{\partial f_j}{\partial z^k}(t,z_*) v^k(t)=0$$
it follows that
$$\frac{d}{d\varepsilon}\Big|_{\varepsilon=0}S(z_*+\varepsilon v)=0.$$
Introduce a linear functional 
$$b:X^{m+n}\to\mathbb{R},\quad b(v)=\frac{d}{d\varepsilon}\Big|_{\varepsilon=0}S(z_*+\varepsilon v),$$
and a linear operator
$$A:X^{m+n}\to X^{l},\quad (Av)_j=\frac{\partial f_j}{\partial z^k}(t,z_*) v^k.$$ It is clear, both these mappings are bounded and 
$$\ker A\subset\ker b.$$
\begin{lemma}\label{xfff}
The operator $A$ maps $X^{m+n}$ onto $X^l$ that is $$A(X^{m+n})=X^l.$$\end{lemma}
{\it Proof.} By $\tilde A(t)$ denote the matrix
$$\frac{\partial f_j}{\partial z^k}\big(t,z_*(t)\big).$$ It is convenient to consider our functions to be defined on the circle $t\in \mathbb S=\mathbb{R}/(\omega\mathbb{Z}).$

Fix an  element $w\in X^l$. Let us cover the circle $\mathbb{S}$ with open intervals $U_i,\quad i=1,\ldots, N$ such that there exists a set of functions $$v_i\in H^1(U_i),\quad \tilde A(t) v_i(t)=w(t),\quad t\in U_i,\quad i=1,\ldots, N.$$
And let $\psi_i$ be a smooth partition of unity subordinated to the covering $\{U_i\}$. A function $\tilde v(t)=\sum_{i=1}^N\psi_i(t) v_i(t)$ belongs to $H^1(\mathbb{S})$ and for each $t$ it follows that $\tilde A(t)\tilde v(t)=w(t)$. But the function $\tilde v$ is not obliged to be odd.

Since $\tilde A(-t)=\tilde A(t)$ we have
$$\tilde A(t) v(t)=w(t),\quad v(t)=\frac{\tilde v(t)-\tilde v(-t)}{2}\in X^{m+n}.$$
The Lemma is proved.

\subsection{}
Recall a lemma from functional analysis \cite{KF}.
\begin{lemma}\label{dfgdfdddd}
Let $E,H,G$ be Banach spaces and $$A:E\to H,\quad B:E\to G$$ be bounded linear operators; $\ker A\subseteq\ker B.$

If the operator $A$ is onto then there exists a bounded operator $\Gamma:H\to G$ such that $B=\Gamma A.$\end{lemma}

Thus there is a linear function $\Gamma\in (X^l)'$ such that $$b(v)=\Gamma A(v),\quad v\in X^{m+n}.$$

Or by virtue of the
Riesz representation theorem, there exists a set of functions $\{\gamma^1,\ldots,\gamma^l\}\subset X$ such that
$$\frac{d}{d\varepsilon}\Big|_{\varepsilon=0}S(z_*+\varepsilon v)=\int_0^\omega\dot\gamma^j(t)\frac{d}{dt}\Big(\frac{\partial f_j}{\partial z^k}(t,z_*)v^k(t)\Big)dt$$
for all $v\in X^{m+n}$.

\subsection{}Every element $v\in X^{m+n}$ is presented as follows 
$$v(t)=\int_0^ty(s)ds,$$ where $y\in Y^{m+n}$ is such that 
$$\int_0^{\omega}y(s)ds=0.$$
Introduce a linear operator $h:Y^{m+n}\to\mathbb{R}^{m+n}$ by the formula
$$h(y)=\int_0^{\omega}y(s)ds.$$
Define a  linear functional $q:Y^{m+n}\to\mathbb{R}$  by the formula
$$q(y)=(b-\Gamma A)v,\quad v(t)=\int_0^t y(s)ds.$$
Now all our observations lead to
$$\ker h\subseteq \ker q.$$
Therefore, there exists a linear functional $\lambda:\mathbb{R}^{m+n}\to\mathbb{R}$ such that
$$q=\lambda h$$

Let us rewrite the last formula explicitly. There are real constants $\lambda_k$ such that for any $y^k\in Y$ one has
\begin{align}
&\int_0^{\omega}\Big(\frac{\partial L}{\partial \dot z^k}(t,z_*,\dot z_*)y^k(t)+\frac{\partial L}{\partial z^k}(t,z_*,\dot z_*)\int_0^ty^k(s)ds\Big)dt\nonumber\\
&=\int_0^\omega\dot\gamma^j(t)\frac{\partial f_j}{\partial z^k}(t,z_*)y^k(t)dt+\int_0^\omega\dot\gamma^j(t)\frac{d}{dt}\Big(\frac{\partial f_j}{\partial z^k}(t,z_*)\Big)\int_0^ty^k(s)dsdt\nonumber\\
&+\lambda_k\int_0^{\omega}y^k(s)ds.\nonumber\end{align}

\subsection{}By the Fubini theorem we obtain
\begin{align}&\int_0^{\omega}\frac{\partial L}{\partial \dot z^k}(t,z_*,\dot z_*)y^k(t)dt+\int_0^{\omega}y^k(s)\int_s^{\omega}\frac{\partial L}{\partial z^k}(t,z_*,\dot z_*)dtds\nonumber\\
&=\int_0^\omega\dot\gamma^j(t)\frac{\partial f_j}{\partial z^k}(t,z_*)y^k(t)dt+\int_0^\omega y^k(s)\int_s^\omega\dot\gamma^j(t)\frac{d}{dt}\Big(\frac{\partial f_j}{\partial z^k}(t,z_*)\Big)dtds\nonumber\\
&+\lambda_k\int_0^{\omega}y^k(s)ds.\label{sfr44}\end{align}
In this formula 
the functions 
$$\frac{\partial L}{\partial \dot z^k}(t,z_*,\dot z_*),\quad \dot\gamma^j(t)\frac{\partial f_j}{\partial z^k}(t,z_*)$$ are even and $\omega$-periodic functions of $t$.

The functions 
$$\dot\gamma^j(t)\frac{d}{dt}\Big(\frac{\partial f_j}{\partial z^k}(t,z_*)\Big),\quad \frac{\partial L}{\partial z^k}(t,z_*,\dot z_*) $$ are odd and $\omega$-periodic.

Employ the following trivial observation.\begin{prop}\label{cfgggtg}If $w\in L^1_{\mathrm{loc}}(\mathbb{R})$ is an $\omega-$periodic and odd function then for any constant $a\in\mathbb{R}$ a function $$t\mapsto\int_a^tw(s)ds$$ is also $\omega-$periodic.\end{prop}
So that the functions 
$$\int_s^\omega\dot\gamma^j(t)\frac{d}{dt}\Big(\frac{\partial f_j}{\partial z^k}(t,z_*)\Big)ds,\quad\int_s^\omega \frac{\partial L}{\partial z^k}(t,z_*,\dot z_*) ds$$
are even and $\omega$-periodic in $s$.

Therefore,  equation (\ref{sfr44}) is rewritten as $(y,\eta)_Y=0$ for any $y=(y^1,\ldots y^{m+n})\in Y^{m+n}$ and $\eta=(\eta_1,\ldots,\eta_{m+n})$ stands for
\begin{align}\eta_k&=\frac{\partial L}{\partial \dot z^k}(t,z_*,\dot z_*)+\int_t^{\omega}\frac{\partial L}{\partial z^k}(s,z_*,\dot z_*)ds\nonumber\\
&-\dot\gamma^j(t)\frac{\partial f_j}{\partial z^k}(t,z_*)-\int_t^\omega\dot\gamma^j(s)\frac{d}{ds}\Big(\frac{\partial f_j}{\partial z^k}(s,z_*)\Big)ds-\lambda_k\nonumber\in Y.\end{align}

Consequently we obtain the following system
\begin{align}&\frac{\partial L}{\partial \dot z^k}(t,z_*,\dot z_*)+\int_t^{\omega}\frac{\partial L}{\partial z^k}(s,z_*,\dot z_*)ds\nonumber\\
&=\dot\gamma^j(t)\frac{\partial f_j}{\partial z^k}(t,z_*)+\int_t^\omega\dot\gamma^j(s)\frac{d}{ds}\Big(\frac{\partial f_j}{\partial z^k}(s,z_*)\Big)ds+\lambda_k.\label{sfddgrf44}\end{align}
 here $ k=1,\ldots,m+n.$ 
 
 If we formally differentiate  both sides of equations (\ref{sfddgrf44}) in $t$ then we obtain the Lagrange equations (\ref{xvvv}) with $\alpha^j=\ddot\gamma^j$.

 Equations  (\ref{sfddgrf44}) hold for almost all $t\in(0,\omega)$ but all the functions contained in (\ref{sfddgrf44}) are defined for all $t\in\mathbb{R}$.

 The functions 
 $$\int_t^\omega\dot\gamma^j(s)\frac{d}{ds}\Big(\frac{\partial f_j}{\partial z^k}(s,z_*)\Big)ds,\quad 
 \int_t^{\omega}\frac{\partial L}{\partial z^k}(s,z_*,\dot z_*)ds$$ 
 are $\omega-$periodic by Proposition \ref{cfgggtg}. Equation (\ref{sfddgrf44}) holds for almost all $t\in\mathbb{R}$.

\subsection{}Let $g^{ij}$ stand for the components of the matrix inverse to $(g_{ij}):\quad g^{ij}g_{ik}=\delta_k^j$. Present equation (\ref{sfddgrf44}) in the form
\begin{align}\label{csr4f}\dot z^j_*(t)&=g^{kj}(t,z_*(t))\nonumber\\
&\cdot\Big(\lambda_k+\dot\gamma^i(t)\frac{\partial f_i}{\partial z^k}(t,z_*(t))+\int_t^\omega\dot\gamma^i(s)\frac{d}{ds}\Big(\frac{\partial f_i}{\partial z^k}(s,z_*(s))\Big)ds\nonumber\\
&-\int_t^{\omega}\frac{\partial L}{\partial z^k}(s,z_*(s),\dot z_*(s))ds -a_k(t,z_*(t))\Big).\end{align}
Together with equation (\ref{csr4f}) consider  equations
\begin{equation}\label{fvvfvfvfsfg}
\frac{\partial f_j}{\partial t}(t,z_*(t))+\frac{\partial f_j}{\partial z^k}(t,z_*(t))\dot z^k_*(t)=0.\end{equation} These equations follow from (\ref{constr}).

Recall that by the Sobolev embedding theorem, $z_*\in X\subset C(\mathbb{R}).$

Due to (\ref{dfggfgfg3e}) we have
$$\det B(t,z_*)\ne 0,\quad B(t,z_*)=\Big(g^{kj}(t,z_*)\frac{\partial f_i}{\partial z^k}(t,z_*)\frac{\partial f_l}{\partial z^j}(t,z_*)\Big)$$ for all $t$. 
Substituting $\dot z_*$ from (\ref{csr4f}) to (\ref{fvvfvfvfsfg}) we can express $\dot\gamma^j$ and see $\dot\gamma^j\in C(\mathbb{R})$. Thus from (\ref{csr4f}) it follows that $\dot z_*\in C(\mathbb{R}).$

Applying this argument again we obtain $\ddot\gamma^j,\ddot z_*\in C(\mathbb{R}).$

This proves the theorem for the case of odd constraints.

\subsection{}
Let us discuss the proof of the theorem under the assumption that the constraints are even (\ref{ncdnvv2}).

\begin{df}
By $Z$ denote a space of functions $u\in   H^1_{\mathrm{loc}}(\mathbb{R})$ such that for all $t\in\mathbb{R}$ the following conditions hold
$$u(-t)=u(t),\quad u(t+\omega)=u(t).$$
\end{df}
The space $Z^l$ is a real Hilbert space with respect to an inner product
$$(u,v)_Z=\sum_{i=1}^l\int_0^\omega\big(u_i(t)v_i(t)+\dot u_i(t)\dot v_i(t)\big)dt.$$ 
This is the standard inner product in $H^1[0,\omega]$.

So what is changed now? The operator $A$ takes the space $X^{m+n}$ onto the space $Z^l$. The proof of this fact is the same as in Lemma \ref{xfff}.

By  the
Riesz representation theorem, there exists a set of functions $\{\gamma^1,\ldots,\gamma^l\}\subset Z$ such that
\begin{align}b(v)&=\int_0^\omega\dot\gamma^j(t)\frac{d}{dt}\Big(\frac{\partial f_j}{\partial z^k}(t,z_*)v^k(t)\Big)dt\nonumber\\
&+\int_0^\omega\gamma^j(t)\frac{\partial f_j}{\partial z^k}(t,z_*)v^k(t)dt\nonumber\end{align}
for all $$v=\int_0^t y(s)ds,\quad y\in Y^{m+n},\quad \int_0^\omega y(s)ds=0.$$
So that equation (\ref{sfddgrf44}) is replaced with the following one
\begin{align}&\frac{\partial L}{\partial \dot z^k}(t,z_*,\dot z_*)+\int_t^{\omega}\frac{\partial L}{\partial z^k}(s,z_*,\dot z_*)ds\nonumber\\
&=\dot\gamma^j(t)\frac{\partial f_j}{\partial z^k}(t,z_*)+\int_t^\omega\dot\gamma^j(s)\frac{d}{ds}\Big(\frac{\partial f_j}{\partial z^k}(s,z_*)\Big)ds\nonumber\\
&+\int_t^\omega\gamma^j(s)\frac{\partial f_j}{\partial z^k}(s,z_*)ds
+\lambda_k.\nonumber\end{align}
 Here $ k=1,\ldots,m+n.$ 
 By the same argument the functions
 
 $$\int_t^\omega\gamma^j(s)\frac{\partial f_j}{\partial z^k}(s,z_*)ds,\quad \int_t^\omega\dot\gamma^j(s)\frac{d}{ds}\Big(\frac{\partial f_j}{\partial z^k}(s,z_*)\Big)ds$$
 are $\omega-$periodic and one can put $\alpha^j=\ddot\gamma^j-\gamma^j.$

 Other argument is the same as above.
The theorem is proved.

\section{Appendix}
\begin{lemma}\label{xzsdffgff}Fix a positive constant $\delta$.

Let $z_1,z_2$ be functions from $ \mathcal V_p$ such that
$$\inf\{|z_i(t)-\sigma'|\mid t\in[0,\omega],\quad \sigma'\in \sigma,\quad i=1,2\}\ge\delta.$$
There exists a positive number $\varepsilon>0$ such that if 
$$\max_{t\in[0,\omega]}|z_1(t)-z_2(t)|<\varepsilon$$ then these functions are homotopic. 
\end{lemma}
{\it Proof of lemma \ref{xzsdffgff}.} Our argument is quite standard. So we present a sketch of the proof.

It is convenient to consider  $z_1,z_2$ as  functions with values in $\mathcal C$ (see Remark \ref{remqq}). In the same sense $F(t)$ is a submanifold in $\mathcal C$ and the functions $z_1,z_2$ define a pair of closed curves in $\mathcal C$. 

Choose  a Riemann metric in $\mathcal C$, for example as follows
$$d\tau^2=\sum_{k=1}^{m+n} (dz^k)^2.$$ This metric inducts a metric  in $F(t)$.

Under the conditions of the Lemma any two points $z_1(t),z_2(t)\in F(t)$ are connected in $F(t)$ with a unique shortest piece of  geodesic
$\chi(\xi,t),\quad \xi\in[0,\tilde \xi(t)]$ such that
$$\chi(0,t)=z_1(t),\quad \chi(\tilde \xi(t),t)=z_2(t).$$
Here $\xi$ is the arc-length parameter.

Define the homotopy as follows $z(s,t)=\chi(s\tilde \xi(t),t),\quad s\in[0,1]$.

The Lemma is proved.

\subsection*{Acknowledgments} The author wishes to thank Professor E. I. Kugushev for useful discussions.

\end{document}